\documentclass[reqno,12pt]{amsart}
\usepackage{a4,latexsym,amssymb,amsfonts,mathrsfs}
\usepackage{amsmath,amscd,verbatim} 
\usepackage{pdfsync} 
\usepackage{enumerate}
\usepackage{amsthm,color}
\usepackage{wasysym}

\usepackage{graphicx,caption,epstopdf,subfigure} 

\usepackage[colorlinks=true, linkcolor=blue]{hyperref}

\usepackage[mathcal]{euscript}
\usepackage{latexsym} 
\usepackage{amsfonts}
\usepackage{amssymb}
\usepackage{amsmath}
\usepackage{braket}
\usepackage{verbatim}
\usepackage{amsthm}
\usepackage[english]{babel}
\usepackage{mathrsfs}
\usepackage{bbm}
\usepackage{amscd} 
\usepackage{enumitem}
\usepackage{todonotes}

\DeclareFontFamily{U}{mathx}{}
\DeclareFontShape{U}{mathx}{m}{n}{<-> mathx10}{}
\DeclareSymbolFont{mathx}{U}{mathx}{m}{n}
\DeclareMathAccent{\widehat}{0}{mathx}{"70}
\DeclareMathAccent{\widecheck}{0}{mathx}{"71}

\thanks{The first author was supported by the Hellenic Foundation for Research and Innovation (H.F.R.I.) under the '2nd Call for H.F.R.I. Research Projects to support Faculty Members \& Researchers' (Project Number: 4662).\\
The second author is a member of Gruppo Nazionale per l'Analisi Matematica, la Probabilità e le loro Applicazioni (GNAMPA) of Istituto Nazionale di Alta Matematica (INdAM)}

\begin{document}
\title[]{\large Embedding Model and de Branges-Rovnyak spaces in Dirichlet Spaces}

\author{Carlo Bellavita}
\email{cbellavita@math.auth.gr}
\address{Department of Mathematics, Aristotle University of Thessaloniki, 54124, Greece }
\author{Eugenio Dellepiane}
\email{eugenio.dellepiane@unimi.it}
\address{Department of Mathematics, Università degli studi di Milano, 20133, Italy}


\begin{abstract}
In this paper we study embeddings between de Branges-Rovnyak spaces $H(b)$ and harmonically weighted Dirichlet spaces $\mathcal{D}(\mu)$ in terms of the boundary spectrum of $b$ and the support of the measure $\mu$, by using elementary reproducing kernel estimates. We completely characterize the embedding between the model spaces $K_u$ and the local Dirichlet spaces $\mathcal{D}_\zeta$, and we discuss some applications. 
\end{abstract}
\maketitle

\thispagestyle{empty}


\ \ \

\setcounter{page}{1}


\theoremstyle{remark}
\newtheorem*{oss}{Remark}

\theoremstyle{remark}
\newtheorem*{example}{Example}

\theoremstyle{remark}
\newtheorem*{question}{Question}

\theoremstyle{plain}
\newtheorem{thm}{Theorem}[section]

\theoremstyle{plain}
\newtheorem{lem}[thm]{Lemma}

\theoremstyle{plain}
\newtheorem{propo}[thm]{Proposition}

\theoremstyle{plain}
\newtheorem{cor}[thm]{Corollary}

\theoremstyle{definition}
\newtheorem{defi}[thm]{Definition}

\newcommand{\N}{\mathbb{N}}
\newcommand{\Z}{\mathbb{Z}}
\newcommand{\Q}{\mathbb{Q}}
\newcommand{\R}{\mathbb{R}}
\newcommand{\C}{\mathbb{C}}
\newcommand{\K}{\mathbb{K}}
\newcommand{\disk}{\mathbb{D}} 

\newcommand{\D}{\mathcal{D}}
\newcommand{\Dm}{\mathcal{D}(\mu)}
\newcommand{\Dn}{\mathcal{D}(\nu)}
\newcommand{\DE}{\mathcal{D}_E}
\newcommand{\Dzf}{\mathcal{D}_{Z(f)}}
\newcommand{\z}{\zeta}
\newcommand{\la}{\lambda}
\newcommand{\om}{\omega}
\newcommand{\Hb}{H(b)}
\newcommand{\Hbarb}{H(\bar{b})}
\newcommand{\M}{M} 
\newcommand{\Lat}{\text{Lat}}
\newcommand{\supp}{\text{supp}}
\newcommand{\Hol}{\text{Hol}}
\newcommand{\sgn}{\text{sgn}}
\newcommand{\HH}{\text{H}}
\newcommand{\LL}{\text{L}}
\newcommand{\de}{\text{ d}}
\newcommand{\af}{\mathfrak{a}}
\newcommand{\T}{\mathbb{T}}
\newcommand{\B}{B}    
\def\CB{\color{black} }
\def\CR{\color{red} } 
\def\CT{\color{cyan} }
\def\CM{\color{magenta} }
\def\CBL{\color{blue} }

\thispagestyle{empty}
\vspace{-1cm}
\section{Introduction}
In this article we deal with spaces of analytic functions on the unit disk $\disk:=\{z\in\C\colon |z|<1\}.$ In this theory, a prominent role is played by the Hardy spaces $H^p(\disk)$, see for example \cite{garnett}. We briefly introduce the different spaces of interest for this work.

\medskip

Given a bounded analytic function $b$ on $\disk$ with $\|b\|_{H^\infty(\disk)}\leq1$, we define the \emph{de Branges-Rovnyak space} $H(b)$ as the reproducing kernel Hilbert space having for reproducing kernel the function
\[ k^b(z,\om):= \frac{1-\overline{b(\om)}b(z)}{1-\overline{\om}z}, \quad \om,z\in\disk.\] 
These spaces were originally introduced by Louis de Branges and James Rovnyak in 1966 as a generalization of the orthogonal complement of the range of multiplication by $b$ on $H^2(\disk)$, see \cite{Debrangesrovnyak}. For a complete introduction to such spaces see \cite{sarHb} and \cite{hb2}.

\medskip

Another space of interest in this paper is the \emph{local Dirichlet space} $\D_\z$. For a fixed point $\z$ on the unit circle $\T:=\partial\disk$, we define the \emph{local Dirichlet integral} at $\z$ of a function $f$ in $\Hol(\disk)$ as
\[ D_\z(f) := \frac{1}{\pi}\int_\disk |f'(z)|^2 \frac{1-|z|^2}{|z-\z|^2} \de A(z),\]
where $\text{d}A$ is the bidimensional Lebesgue measure. We call $\D_\z$ the space of functions $f$ in $\Hol(\disk)$ such that $D_\z(f)<\infty$. These spaces are studied in \cite{RS} and they belong to a more general class, the so-called \emph{harmonically weighted Dirichlet spaces}. We will discuss this later.

\medskip

In Section $3$ we provide a sufficient condition and a necessary one in order to have an embedding between de Branges-Rovnyak spaces and local Dirichlet spaces, i.e. a bounded inclusion $H(b)\hookrightarrow\D_\z$. Both these conditions involve the notion of \emph{boundary spectrum}: given a bounded analytic function $b$ with $\|b\|_{H^\infty(\disk)}=1,$ we define its boundary spectrum as the set
\[\sigma(b):= \{\la\in\T\colon \liminf_{z\to\la} |b(z)|<1\}.\] 
As we will explain later, this set carries information about regularity of the function $b$ and of all elements of $H(b)$. In particular, we proved the following results.
\begin{thm}
\label{HbinDz}
Let $b$ be a bounded analytic function with $\|b\|_{H^\infty(\disk)}=1$, and let $\z\in\T$ be such that $\z\notin\overline{\sigma(b)}$. Then, the embedding $H(b)\hookrightarrow \D_\z$ holds.
\end{thm}

\begin{thm}\label{KuinDzb}
Let $b$ be a bounded analytic function with $\|b\|_{H^\infty(\disk)}=1$, and let $\z\in\T$ be such that $\z\in\sigma(b)$. Then, the de Branges-Rovnyak space $H(b)$ does not embed into the local Dirichlet space $\D_\z$.
\end{thm}

Later in the article, we restrict our attention to the \emph{model spaces}. Given an inner function $u$, i.e. a bounded analytic function on $\disk$ with $|u|=1$ a.e. on $\T$, we define the model space $K_u$ as the complementary space $K_u := H^2(\disk)\ominus uH^2(\disk)$. These spaces naturally arise as the closed invariant subspaces of the backward shift operator $S^*$ on $H^2(\disk)$. For a complete introduction on this subject, we refer to \cite{nikolski} and \cite{1}. The model spaces are a particular class of de Branges-Rovnyak spaces: when one considers an inner function $u$, one has that $H(u)=K_u$ with equality of norms. 

\medskip

As we mentioned before, we also consider the class of harmonically weighted Dirichlet spaces $\Dm$. Given a finite positive Borel measure $\mu$ on the unit circle $\T$, the associated $\Dm$ space is the space of holomorphic functions on $\disk$ having finite \emph{harmonically weighted} Dirichlet integral
\begin{equation}\label{Dm}
    D_\mu (f) :=\frac{1}{\pi} \int_\disk |f'(z)|^2 \, P\mu(z) \de A(z),
\end{equation} 
where $P\mu$ is the Poisson integral of $\mu$,
\[ P\mu(z) := \int_\T \frac{1-|z|^2}{|\lambda -z|^2} \de\mu(\lambda), \qquad z \in \disk.\]
These spaces were introduced by Stefan Richter in 1991 for the representation of cyclic analytic two-isometries, see \cite{R}. Also, they play a key role in the description of the closed shift-invariant subspaces of the classical Dirichlet space $\D := \D(m)$, where $m$ is the Lebesgue measure on $\T$, see \cite{RS}. In Section $5$ we deal with the embedding $K_u \hookrightarrow \Dm$. Again, we provide a sufficient condition and a necessary one for the embedding to hold, involving the support of the measure $\mu$ and the boundary spectrum $\sigma(u)$ of the inner function $u$.
\begin{thm}
\label{KuinDm}
Let $\mu$ be a finite positive Borel measure on $\T$ and let $u$ be an inner function. If $\supp(\mu)\cap\sigma(u) = \emptyset$, then the embedding $K_u \hookrightarrow \Dm$ holds.
\end{thm}

\begin{thm}\label{KuinDmnec}
Let $\mu$ be a finite positive Borel measure on $\T$ and let $u$ be an inner function. If the embedding $K_u \hookrightarrow \Dm$ holds, then $\mu\big(\sigma(u)\big)=0$.
\end{thm}

The paper is organized as follows. Section $2$ is devoted to some well-known preliminaries. In Section $3$ we describe the embedding $H(b)\hookrightarrow \D_\zeta$ for general $b$'s. In Section $4$ we discuss some applications of the embedding $K_u\hookrightarrow\D_\z$. In the fifth section, we prove Theorems \ref{KuinDm} and \ref{KuinDmnec}. We conclude with an open problem.

\section{Preliminaries}
We introduce the main spaces involved in this article. Let us start with the harmonically weighted Dirichlet spaces. Given a finite positive Borel measure $\mu$ on the unit circle $\T$, its Poisson integral is the harmonic function 
\[ P\mu(z) := \int_\T \frac{1-|z|^2}{|\lambda -z|^2} \de\mu(\lambda), \qquad z \in \disk.\]
The associated harmonically weighted Dirichlet space $\Dm$ is   
\[ \Dm := \{ f\in\Hol(\disk) \colon D_\mu(f) <\infty\},\]
where 
\begin{equation}\label{Dm}
    D_\mu (f) :=\frac{1}{\pi} \int_\disk |f'(z)|^2 \, P\mu(z) \de A(z)
\end{equation} 
is the {harmonically weighted Dirichlet integral}. Notice that $D_\mu$ is a seminorm that annihilates the constants. We recall a few basic properties; for a treatise of Dirichlet spaces we refer to \cite{primer}. If $\mu$ is a finite measure on $\T$ such that $\mu(\T) >0$, then $\Dm$ is a subset of $H^2(\disk)$ which contains all polynomials. Moreover, $\Dm$ is a Hilbert space with respect to the inner product induced by the norm 
\[ \|f\|_\mu^2 := \|f\|_{H^2}^2 + D_\mu(f).\]
For $\z\in\T$, considering the Dirac delta $\delta_\z$ we obtain the so-called local Dirichlet space, which we simply denote  by $\D_\z$. Also, we write $D_\z(f)$ instead of $D_{\delta_\z}(f)$. For $f \in H^2(\disk)$, by Fubini’s theorem, $D_\mu(f)$ given in \eqref{Dm} can be expressed as 
\begin{equation}\label{eqfubini}
    D_\mu(f) = \int_{\T} D_\zeta(f) \, d\mu(\zeta).
\end{equation} 
In  \cite{RS} Richter and Sundberg proved the following useful formula for $D_\z(f)$, which includes the boundary value $f(\z)$ defined as the radial limit $\lim_{r\to1^-}f(r\z)$, whenever it exists.
\begin{thm}[Local Douglas formula]
\label{localdouglasthm}
Let $f\in H^2(\disk)$ and $\z\in\T$. If the boundary value $f(\z)$ exists, then 
\begin{equation}
\label{localdouglas}
D_\z(f)=\int_\T \,\,\bigg| \frac{f(\la)-f(\z)}{\la -\z}\bigg|^2 \de m(\la).
\end{equation}
On the other hand, if $f(\z)$ does not exist, then $D_\z(f)=\infty$. In particular, all functions in $\D_\z$ admit boundary value at $\z$.
\end{thm} 
This formula shows that the quotient ratio at $\z$ plays an important role for membership in the local Dirichlet space. Richter and Sundberg also proved the following characterization of $\D_\zeta$. One has
\begin{equation}\label{second defn}
    \D_\zeta =\left\lbrace f \in \text{Hol}(\disk)\ : f(z)=c+(z-\z)g(z), \text{where } c \in \mathbb{C} \text{ and } g\in H^2(\disk)\right\rbrace ,
\end{equation}
with the equality $D_\z(f)=\| g\|_{H^2}^2$. Some aspects of the structure of local Dirichlet spaces have been recently studied by Fricain and Mashreghi in \cite{7}. 
Finally, we point out that local Dirichlet spaces of order $m\in\N$ have been recently introduced by Luo, Gu and Richter in \cite{LGR} and further developed in \cite{9} and \cite{8}. 

\medskip

In the rest of this section we provide some preliminary information about de Branges-Rovnyak spaces. There are many equivalent ways to define these spaces: we will follow the reproducing kernel approach. As shown in the classic work of Aronszajn in \cite{aronszajn}, given a positive definite function $k$ on $\disk\times\disk$ one can construct a Hilbert space $H_k$ of functions on $\disk$ such that for all $\om\in\disk$ the function $k(\cdot,\om)$ belongs to $H_k$ and it holds the so-called \emph{reproducing kernel property}, i.e. 
\[ f(\om) = \langle f,k(\cdot,\om) \rangle_{H_k}, \quad f\in H_k.\]
Given a bounded analytic function $b$ on $\disk$ with $\|b\|_{H^\infty(\disk)}\leq1$, the de Branges-Rovnyak space $H(b)$ is the reproducing kernel Hilbert space having for reproducing kernel the function
\[ k^b(z,\om):= \frac{1-\overline{b(\om)}b(z)}{1-\overline{\om}z}, \quad \om,z\in\disk.\]

We denote by $\langle \cdot,\cdot\rangle_b$ the inner product of $H(b)$ and by $\|\cdot\|_b$ its induced norm.
$H(b)$ is a space of analytic functions contained in $H^2(\disk)$ and it holds the norm inequality
\begin{equation}
\label{HbinH2}
    \|f\|_{H^2} \leq \|f\|_{b}, \qquad f\in H(b).
\end{equation}
However, in general $H(b)$ is not complete with respect to the $H^2$ norm. If $b=u$ is an inner function, then $H(u)$ coincides with the \emph{model space} $K_u$, defined as the orthogonal complement $K_u:=H^2(\disk) \ominus uH^2(\disk)$. Therefore, $H(u)=K_u$ is closed in $H^2(\disk)$, and it holds the norm identity $\|\cdot\|_b=\|\cdot\|_{H^2}$. As a corollary of a classic result of Beurling (see Theorem $4.3$ in \cite{1}), the closed $S^\ast$-invariant subspaces of $H^2(\disk)$ are exactly the model spaces. More in general, all de Branges-Rovnyak spaces are $S^\ast$-invariant. The operator 
\[X_b : H(b) \ni f \mapsto S^\ast f \in H(b)\]
is well-defined and bounded.

\medskip 

In order to introduce the notion of boundary spectrum, first we recall a key factorization result (see for example Theorem $3.20$ in \cite{1}): every analytic function $b$ with $\|b\|_{H^\infty(\disk)} = 1$ can be factorized as $b=Ou$, where $O$ is the outer function
\begin{equation}\label{hinf function}
    O(z):= \exp{\bigg\{\int_\T \frac{\z+z}{\z-z} \log|b(\z)| \de m(\z)\bigg\}}
\end{equation}
and $u$ is an inner function. In particular, according to the Nevanlinna factorization, we can write
\begin{equation}\label{inner function}
    u(z) = \prod_{n=1}^\infty \frac{\overline{a_n}}{|a_n|}\frac{a_n-z}{1-\overline{a_n}z} \,\exp{\bigg\{-\int_\T \frac{\z+z}{\z-z} \de \tau(\z)\bigg\}},
\end{equation}
where $\{a_n\}_{n\geq 1}$ are the zeros of $u$ and $\tau$ a positive singular measure.

\begin{defi}
For a bounded analytic function $b$ on $\disk$ with $\|b\|_{H^\infty(\disk)} = 1$, we define its {boundary spectrum} as the set
\[ 
\sigma(b) := \{\la\in\T \colon \liminf_{z\to\la} |b(z)|<1\}.
\]
\end{defi}
As stated in \cite{10}, the closure $\overline{\sigma(b)}$ is the smallest closed subset of $\T$ containing the closure of the zero set $\{a_n\}_n$ and the supports of the (positive finite) measures $\tau$ and $-\log|b(\zeta)|dm(\zeta)$. It is known that $b$ has an analytic extension through any arc of the open set $\T \setminus \overline{\sigma(b)}$ with unimodular values on such arcs, see again \cite{10}. 
If $b=u$ is an inner function, then it holds
\[ \sigma(u)=\{\la\in\T \colon \liminf_{z\to\la} |u(z)|=0\}.\]
In particular, the spectrum of inner functions is a closed set. We also note that there exist bounded functions with closed spectrum that are not necessarily inner, for example one-component bounded functions, defined and studied in \cite{12}. The name \emph{spectrum} comes from the following fact. 

\begin{thm}\label{spectrumshift}
Let $b$ be a bounded analytic function on $\disk$ with $\|b\|_{H^\infty(\disk)}=1$. Then, the intersection of the spectrum of the operator $X_b^*$ and the unit circle $\T$ coincides with the closure of the boundary spectrum of $b$. In symbols, 
\[ \sigma(X_b^*) \cap \T = \overline{\sigma(b)}. \]
\end{thm}

 For a proof, see Corollary $20.14$ in \cite{hb2}. However, we point out that the definition of $\sigma(b)$ used in \cite{hb2}, found in the first volume of the same book \cite{hb1}, is different from the one used in this paper, taken from \cite{10}. In particular, in this paper $\sigma(b)$ is not necessarily closed.

The boundary regularity of the function $b$ results in properties of functions in $H(b)$. The notion we need is the \emph{angular derivative in the sense of Caratheodory} (ADC). We say that an analytic function $b$ on $\disk$ with $\|b\|_{H^\infty(\disk)}\leq1$ admits ADC at $\z\in\T$ if the derivative $b'$ admits non-tangential limit at $\z$ and $|b(\z)|=1$. The result that follows is Theorem $21.1$ in \cite{hb2}.

\begin{thm}
\label{thm211debranges}
Let $b$ be an analytic function on $\disk$ with $\|b\|_{H^\infty(\disk)}\leq1$ and let $\z \in \T$. 
The following are equivalent:
\begin{enumerate}[label=(\roman*)]
\item There exists $\la \in\T$ such that the function 
\[ \disk \ni z \mapsto \frac{b(z)-\la}{z-\z}\]
belongs to $H(b)$.
\item Every function $f$ in $H(b)$ admits non-tangential limit at $\z$.
\item $b$ has ADC at $\z$.
\end{enumerate}
Furthermore, under these conditions, $\la=b(\z)$ and for every $f\in H(b)$ one has $f(\z) = \langle f,k_{\z}^b\rangle_b$, where
\[ k_\z^b(z) = \frac{1-\overline{b(\z)}b(z)}{1-\overline{\z}z}\in H(b).\]
\end{thm}

Also, by Theorem $18.21$ in \cite{hb2}, the operator $X^*_b$ intertwines the reproducing kernels, in the sense that
\[ k_z^b = (I-\overline{z} X^*_b)^{-1}k_0^b, \qquad z\in\disk.\]
One can easily prove that the same formula still holds replacing $z\in\disk$ with $\z\in\T\setminus\overline{\sigma(b)}$, that is
\begin{equation}\label{rel kernel}
k_\z^b = (I-\overline{\z} X^*_b)^{-1}k_0^b, \qquad \z\in\T\setminus\overline{\sigma(b)}.
\end{equation}

\section{The embedding $H(b) \hookrightarrow \D_\zeta$}
We come now to our two main results.

\proof[Proof of Theorem \ref{HbinDz}]

By assumption, $b$ extends analytically in a neighbourhood of $\z$. Then, by Theorem \ref{thm211debranges}, every function in $H(b)$ admits non-tangential boundary value at $\z$. Also, since $\z\notin \overline{\sigma(b)}$, by Theorem \ref{spectrumshift} the operator $I-\overline{\z} X^*_b$ is boundedly invertible in $H(b)$ and, by \eqref{rel kernel}, 
\[ 
k_\z^b = (I-\overline{\z} X^*_b)^{-1}k_0^b. 
\]
Thus, the operator
\begin{equation}\label{Q_u}
    Q_\z^b := ( I-\z X_b)^{-1} X_b 
\end{equation}
is bounded on $H(b)$. By an algebraic computation, for $z\in\disk$ it holds the operator identity
\[
(I-zX_b)^{-1}(I-\z X_b)^{-1} X_b=\frac{(I-zX_b)^{-1}-(I-\z X_b)^{-1}}{z-\z}.
\]
For every $f\in H(b)$ it holds the formula
\[ 
Q_\z^b f(z) = \frac{f(z)-f(\z)}{z-\z}, \qquad z\in\disk.
\] 
This formula for the operator $Q_\omega^b$, with $\omega\in\mathbb{D}$, is found in \cite[Chapter $2$]{sarHb}. It continues to hold for $\zeta\notin\overline{\sigma(b)}$ since
\begin{align*}
Q_\z^b f (z) &= \langle Q_\z^b f,k_z^b \rangle_b  = \langle Q_\z^b f, (I-\overline{z} X^*_b)^{-1}k_0^b \rangle_b =  \langle (I-zX_b)^{-1} Q^b_\z f, k_0^b\rangle_b \\
&= \langle (I-zX_b)^{-1}(I-\z X_b)^{-1} X_b f,k_0^b\rangle_b = \frac{1}{z-\z} \langle (I-zX_b)^{-1}f - (I-\z X_b)^{-1} f,k_0^b\rangle_b \\
&= \frac{1}{z-\z} \langle f, (I-\overline{z}X^*_b)^{-1} k_0^b\rangle_b - \frac{1}{z-\z} \langle  f,(I-\overline{\z}X^*_b)^{-1}k_0^b\rangle_b \\
&= \frac{1}{z-\z} \langle f, k_z^b \rangle_b - \frac{1}{z-\z} \langle  f, k_\z^b \rangle_b = \frac{f(z)-f(\z)}{z-\z}.
\end{align*}
This proves the boundedness of the embedding $H(b) \hookrightarrow \D_\z$, since by Theorem \ref{localdouglasthm} and \eqref{HbinH2}
\[ \|f\|_{H^2}^2 + D_\z(f) = \|f\|_{H^2}^2 + \|Q_\z^b f\|_{H^2}^2 \leq  \|f\|_b^2 + \|Q_\z^b f\|_b^2 \leq \big( 1+ \|Q_\z^b\|^2\big) \|f\|_b^2. \qed
\]

\medskip 

\begin{proof}[Proof of Theorem \ref{KuinDzb}]
By contradiction, let us suppose that the embedding $H(b) \hookrightarrow \D_\z$ holds. Let $C>0$ be such that
\begin{equation}
\label{boundassurdob}
 D_\z(f) \leq C\|f\|_b^2, \qquad f\in H(b).
\end{equation}
By assumption, $\z\in\sigma(b)$, hence there exists a sequence $(\om_n)_n$ in $\disk$ converging to $\z$ such that 
\[ \beta:= \lim_n |b(\om_n)|<1.\]
Let us consider the family of kernels 
\[ k_n(z) := k_{\om_n}^b(z) = \frac{1-\overline{b(\om_n)}b(z)}{1-\overline{\om_n}z}.\]
Since $H(b) \subseteq \D_\z$, by Theorem \ref{localdouglasthm} every function of $H(b)$ admits boundary value at $\zeta$. By Theorem \ref{thm211debranges}, $b(\zeta)$ is well defined and unimodular. Therefore, one can compute
\begin{align*}
k_n(z) - k_n(\z) &= \frac{1-\overline{b(\om_n)}b(z)}{1-\overline{\om_n}z} - \frac{1-\overline{b(\om_n)}b(\z)}{1-\overline{\om_n}\z} \\
&= \frac{\overline{\om_n}(z-\z) - \overline{b(\om_n)}\big(b(z)-b(\z)\big) + \overline{\om_n b(\om_n)}\big(\zeta b(z)-zb(\z)\big)}{(1-\overline{\om_n}z)(1-\overline{\om_n}\z)}\\
&=\frac{\overline{\om_n}(z-\z)(1-\overline{b(\omega_n)}b(\z)) - \overline{b(\om_n)}\big(b(z)-b(\z)\big)(1-\overline{\omega_n}\z) }{(1-\overline{\om_n}z)(1-\overline{\om_n}\z)}
.
\end{align*}
Consequently,
\begin{align}
    \nonumber \frac{k_n(z) - k_n(\z)}{z-\z} &=
\frac{\overline{\om_n}}{1-\overline{\om_n}z}\frac{1-\overline{b(\om_n)}b(\z)}{1-\overline{\om_n}\z} - \frac{\overline{b(\om_n)}}{1-\overline{\om_n}z}\frac{b(z)-b(\z)}{z-\z}\\
\label{quotientku2} &= \overline{\om_n} c_{\om_n}(z) k_n(\z)   - \overline{b(\om_n)} c_{\om_n}(z) b(\z) \overline{\z}  k_\z^b(z),
\end{align}
where  \[ c_{\om_n}(z) = \frac{1}{1-\overline{\om_n}z}\] is  the usual Szeg\"o kernel, the reproducing kernel of the Hardy space $H^2(\disk)$. The local Dirichlet integral can be computed as in \eqref{localdouglas}, yielding
\begin{align*}
 D_\z(k_n) &= \bigg\| \frac{k_n-k_n(\z)}{\cdot-\z}\bigg\|_{H^2}^2 \\
&= \big\langle \overline{\om_n}\, k_n(\z) c_{\om_n}  - \overline{b(\om_n)} b(\z) \overline{\z} c_{\om_n} k_\z^b, \overline{\om_n}\, k_n(\z) c_{\om_n}  - \overline{b(\om_n)} b(\z) \overline{\z} c_{\om_n} k_\z^b \big\rangle_{H^2} \\
&= |\om_n|^2 |k_n(\z)|^2 \|c_{\om_n}\|_{H^2}^2 -2\Re\Big(\overline{\om_n}k_n(\z) b(\om_n) \overline{b(\z)}\z \langle c_{\om_n},c_{\om_n}k_\z^b\rangle_{H^2} \Big) +|b(\om_n)|^2 \|c_{\om_n}k_\z^b\|_{H^2}^2. 
\end{align*}
We have written the local Dirichlet integral $D_\z(k_n)$ as a sum of three terms. We leave the first one as it is and work on the other two. 
We use the reproducing property of the Szeg\"o kernel, the fact that $c_{\om_n}k_\z^b$ is an $H^2$ function and we estimate the real part with the modulus, obtaining
\begin{align*}
\Re\Big(\overline{\om_n}k_n(\z) b(\om_n) \overline{b(\z)}\z \langle c_{\om_n},c_{\om_n}k_\z^b\rangle_{H^2} \Big) &= \Re\Big(\overline{\om_n}k_n(\z) b(\om_n) \overline{b(\z)}\z \overline{c_{\om_n}(\om_n)k_\z^b(\om_n) }\Big) \\
&=\|c_{\om_n}\|_{H^2}^2 \Re\big(\overline{\om_n} b(\om_n) \overline{b(\z)}\z k_n(\z)^2 \big) \\
&\leq \|c_{\om_n}\|_{H^2}^2 |k_n(\z)|^2 \, |\om_n b(\om_n)|.
\end{align*}
For the third summand, using the triangular inequality, we have
\begin{align*}
\|c_{\om_n}k_\z^b\|_{H^2}^2 &= \int_\T \bigg| \frac{1}{1-\overline{\om_n}\la}\frac{1-\overline{b(\z)}b(\la)}{1-\overline{\z}\la}\bigg|^2 \de m(\la) \\
&=  \int_\T \bigg|\frac{1}{1-\om_n\overline{\la}} \frac{1}{1-\overline{\om_n}\la}\bigg(\frac{1-\overline{b(\z)}b(\la)}{1-\overline{\z}\la}\bigg)^2\bigg| \de m(\la) \\
&\geq  \bigg| \int_\T \frac{1}{1-\om_n\overline{\la}} \frac{1}{1-\overline{\om_n}\la}\bigg(\frac{1-\overline{b(\z)}b(\la)}{1-\overline{\z}\la}\bigg)^2 \de m(\la) \bigg|.
\end{align*}
The function $c_{\om_n} \big(k_\z^b\big)^2$ belongs to $H^1(\disk)$ and in particular the Cauchy integral formula holds
\begin{align*}
 \int_\T \frac{1}{1-\om_n\overline{\la}} \frac{1}{1-\overline{\om_n}\la}\bigg(\frac{1-\overline{b(\z)}b(\la)}{1-\overline{\z}\la}\bigg)^2 \de m(\la) &= \int_\T \frac{c_{\om_n}(\la) \big(k_\z^b\big)^2(\la)}{1-\om_n\overline{\la}} \de m(\la) \\ 
&= c_{\om_n}(\om_n) \big(k_\z^b\big)^2(\om_n).
\end{align*}
Using this, we obtain
 \[ |b(\om_n)|^2 \|c_{\om_n}k_\z^b\|_{H^2}^2 \geq |b(\om_n)|^2 \|c_{\om_n}\|_{H^2}^2 |k_n(\z)|^2.\]
Now, computing the norms of the kernels 
\[ \|c_{\om_n}\|_{H^2}^2 = \frac{1}{1-|\om_n|^2}, \qquad \|k_n\|_b^2 = \frac{1-|b(\om_n)|^2}{1-|\om_n|^2},\]
we obtain the lower bound 
\begin{align*}
 \frac{D_\z(k_n)}{\|k_n\|_b^2} &\geq \frac{|\om_n|^2 |k_n(\z)|^2}{1-|b(\om_n)|^2} -\frac{2 |k_n(\z)|^2|\om_n b(\om_n)|}{1-|b(\om_n)|^2} + \frac{|b(\om_n)|^2 |k_n(\z)|^2}{1-|b(\om_n)|^2}\\
&= |k_n(\z)|^2\frac{|\om_n|^2 -2 |\om_n b(\om_n)| + |b(\om_n)|^2}{1-|b(\om_n)|^2} \\
&=\bigg| \frac{1-\overline{b(\om_n)}b(\z)}{1-\overline{\om_n}\z}\bigg|^2\frac{\big( |\om_n| - |b(\om_n)| \big)^2}{1-|b(\om_n)|^2} \\
&\geq  \frac{\big(1-|b(\om_n)|\big)^2}{|1-\overline{\om_n}\z|^2} \frac{\big( |\om_n| - |b(\om_n)| \big)^2}{1-|b(\om_n)|^2}.
\end{align*}
Since $\lim_n \om_n = \z$ and $\lim_n |b(\om_n)|=\beta \in [0,1),$ we conclude that
\[ \liminf_n \frac{D_\z(k_n)}{\|k_n\|_b^2} \geq \liminf_n \frac{(1-\beta)^2}{|1-\overline{\om_n}\z|^2} \frac{(1-\beta)^2}{1-\beta^2}= +\infty  ,\]
contradicting the uniform bound in \eqref{boundassurdob}.
\end{proof}
\begin{oss}
In Theorem \ref{KuinDzb} it is shown that there cannot be a (bounded) embedding $H(b)\hookrightarrow \D_\z$, if $\z\in\sigma(b)$. By the closed graph theorem, even a set inclusion $H(b) \subseteq \D_\z$ cannot hold.
\end{oss}

The following result is contained in the proof of Theorem \ref{KuinDzb}.

\begin{cor}
Let $b$ be an analytic function on $\disk$ with $\|b\|_{H^\infty(\disk)}\leq1$ and let $\z \in \T$. If $b$ admits ADC at a point $\z\in\T$, then for all $\om\in\disk$ the reproducing kernel $k^b_\omega$ belongs to $\D_\zeta$.
\end{cor} 
\begin{proof}
From \eqref{quotientku2}, it follows that
\[D_\zeta(k^b_\omega) =    
     \bigg\| \frac{k^b_\omega-k^b_\omega(\z)}{\cdot-\z}\bigg\|_{H^2}^2=\|\overline{\omega}  k_\omega^b(\z) c_{\om}  - \overline{b(\omega)}  b(\z) \overline{\z} c_{\omega}  k_\z^b\|_{H^2}^2 <\infty.\]
\end{proof}

  
We have proved a positive result, that is, that $H(b) \hookrightarrow \D_\z$ when $\z\notin\overline{\sigma(b)}$, and a negative one, that is, that if $\z\in\sigma(b)$, $H(b) \nsubseteq \D_\z$. We now present some examples to show that for the remaining case $\z\in\overline{\sigma(b)}\setminus\sigma(b)$, anything can happen.

\begin{example}
Set
\begin{equation}
    \label{constantw0}
     w_0 := \frac{3-\sqrt{5}}{2},
\end{equation}
let $\z\in\T$ and define the function
\[ b_\z(z) = \frac{(1-w_0)\overline{\z} z}{1-w_0\overline{\z} z}.\]
By Proposition 2 in \cite{5}, $H(b_\z) = \D_\z$ with equality of norms, guaranteeing the embedding. Since $b_\z$ is continuous up to the boundary $\T$, it holds 
\[\sigma(b_\z) =\{ \la\in\T\colon |b_\z(\la)|<1\}. \]
Writing $\z=e^{i\eta}$ and $\la=e^{i\theta}$, one can easily see that
\[ |1-w_0 \overline{\z}\la|^2 = 1-2w_0\cos(\theta-\eta)+w_0^2 > 1-2w_0+w_0^2 = |(1-w_0) \overline{\z}\la|^2, \qquad \text{if} \quad e^{ i \theta}\neq e^{i\eta}, \]
whereas
\[ |1-w_0 \overline{\z}\la|^2 =  |(1-w_0) \overline{\z}\la|^2, \qquad \text{if} \quad e^{ i \theta}= e^{i\eta}.\]
This means that $\sigma(b_\z) = \T \setminus \{\z\}$,   providing a function in $H^\infty(\disk)$ such that $\z\in\overline{\sigma(b_\z)}\setminus\sigma(b_\z)$ while the embedding $H(b_\z)\hookrightarrow\D_\z$ holds.
\end{example}
 
Now we provide an example of a case with $1\in\overline{\sigma(b)}\setminus\sigma(b)$ such that $H(b)\hookrightarrow \D_1$ doesn't hold. We use the following proposition as a criterion for the inclusion, see Corollary $27.18$ in \cite{hb2}:
\begin{propo}
\label{cormashreghi}
    Suppose $b_1$ is a non-extreme point of the closed unit ball of $H^\infty(\disk)$, and assume $b_1$ is continuous on the closed unit disk. Let $b_2$ be a function in $H^\infty(\disk)$ and $\theta_2$ its inner factor. Then the following are equivalent:
    \begin{enumerate}[label=(\roman*)]
        \item It holds the inclusion of de Branges-Rovnyak spaces $H(b_2) \subset H(b_1)$.
        \item The following conditions hold:
        \begin{itemize}
            \item $\{\la\in\T\colon |b_1(\la)| =1\}\cap \sigma(\theta_2) = \emptyset$.
            \item There exists $\gamma >0$ such that $1-|b_2|^2 \leq \gamma (1-|b_1|^2)$ a.e. on $\T$.
        \end{itemize}
    \end{enumerate}
\end{propo}

\begin{example}
    Let \[b_1(z):=  \frac{(1-w_0) z}{1-w_0 z},\] where $w_0$ is the constant in \eqref{constantw0}, so that $H(b_1)=\D_1$, and we construct an outer function $b_2$ as follows. We start by considering the function $\varphi$ defined on $\T$ as
    \[ \varphi(\la)= 
    \begin{cases} \log\Big(\sqrt{1-|1-\la|^{\frac{3}{2}}}\Big), & \text{if} \quad |\arg(\la)|\leq\frac{\pi}{6} , \\
    0, & \text{elsewhere}.
    \end{cases}\]
     The function $\varphi$ is in $L^{\infty}(\T)$ and real-valued, and this allows us to define the outer function    
    \[ b_2(z) := \exp\bigg\{ \int_\T \frac{\la+z}{\la-z} \, \varphi(\la) \!\de m(\la) \bigg\},\]
that satisfies $|b_2| = e^{\varphi}$ a.e. on $\mathbb{T}$. The first condition of $(ii)$ in Proposition \ref{cormashreghi} is trivially true, since $b_2$ is outer and therefore $\sigma(\theta_2)=\emptyset$. For the second condition of $(ii)$, it holds that 
\[  1-|b_2(\la)|^2= |1-\la|^{\frac{3}{2}},\qquad \text{for a.e.}\,\, \lambda\in\mathbb{T}\,\,\text{with} \,\,|\arg(\la)|<\frac{\pi}{6}. 
    \] 
    Since for all $\la\in\T$ it holds
    \[ 1-|b_1(\la)|^2 = \frac{(1-w_0)^2|1-\la|^2}{|1-w_0\la|^2},\]
    it follows that in proximity of the point $1$ the condition $(ii)$ of Proposition \ref{cormashreghi} fails, meaning that the inclusion $H(b_2) \subset H(b_1) = \D_1$ cannot hold. Finally, from a classical argument with Poisson kernels found in the proof of \cite[Theorem $1.9$]{1}, it follows that for every $\lambda\in\mathbb{T}$ with $|\arg{\lambda}|<\frac{\pi}{6}$, it holds that  \[\lim_{z\to\lambda} |b_2(z)| = e^{\varphi(\lambda)} = \sqrt{1-|1-\lambda|^{\frac{3}{2}}},\]
    since $\varphi$ is continuous on such $\lambda$'s and bounded on $\mathbb{T}$. It follows that
    \[\sigma(b)\cap\{\la\in\T\colon|\arg(\la)|<\pi/6\}= \{\la\in\T\colon|\arg(\la)|<\pi/6\}\setminus\{1\}, \] 
    so that $1\in\overline{\sigma(b_2)}\setminus\sigma(b_2)$ while $H(b_2) \not\subset \D_1$.
\end{example}
 
\vspace{10 pt}
\section{Applications of $K_u \hookrightarrow \D_\zeta$}
In the last two sections, we focus on the model space $K_u$. 

\begin{cor}
\label{KuinDz}
    Let $b$ be an analytic function with $\|b\|_{H^\infty}=1$ with closed boundary spectrum, and let $\z\in\T$. Then, the embedding $H(b)\hookrightarrow\D_\z$ holds if and only if $\z\notin\sigma(b)$. In particular, if $u$ is an inner function, then the embedding $K_u\hookrightarrow\D_\z$ holds if and only if $\z\notin\sigma(u)$. 
\end{cor}

\begin{proof}
The result follows using Theorems \ref{HbinDz} and \ref{KuinDzb} and the fact that the spectrum of $b$ is closed.
\end{proof}

We can rewrite the embedding $K_u \hookrightarrow \D_\z$ in terms of the boundedness of the derivative operator, providing a corollary which is somehow related to the results of Baranov about the boundedness of the differentiation operator acting on model spaces, see \cite{BARANOV2013541}.
\begin{cor}
\label{CorDer} Let $u$ be an inner function and $\z\in\T$. Let $D$ be the derivative operator 
\[ D \colon K_u \to  L^2\big(P\delta_\z \!\de A \big), \quad f \mapsto f',\]
acting from the model space to the Lebesgue space $L^2(\disk, P\delta_\z \!\de A)$. Then, $D$ is bounded if and only if $\z\notin\sigma(u)$.
\end{cor}

\proof
It follows at once from Corollary \ref{KuinDz}, for 
\[ \|f'\|_{L^2\big(P\delta_\z \!\de A \big)} =\int_\disk |f'(z)|^2 \frac{1-|z|^2}{|z-\z|^2} \de A(z) = \pi D_\z(f). \qed\]

\medskip

As already said in the introduction, the embedding $K_u \hookrightarrow \D_\zeta$ allows one to find some Carleson measures for $K_u$. First, let us recall the definition.
\begin{defi}
    Let $H$ be a Hilbert space of holomorphic functions on $\disk$. We say that a positive Borel measure $\nu$ on $\disk$ is a \emph{Carleson measure} for $H$ if there exists a constant $C>0$ such that 
\begin{equation}
\label{Carlesondef}
    \int_\disk |f|^2 \de\nu\leq C\|f\|_H^2, \qquad f\in H.
\end{equation}
\end{defi}
Carleson measures for $H^2(\disk)$ appeared in a very natural and powerful way in the proof of the Corona Theorem for $H^\infty(\disk)$, see \cite{garnett}. Such measures have been well studied, and they admit a nice geometric characterization in terms of \emph{Carleson boxes}.

\begin{propo}
\label{carmeasH2}
    Let $\nu$ be a finite positive Borel measure on $\disk$. Given an arc $I\subseteq\T$, the Carleson box associated to $I$ is 
    \[ S(I) := \{re^{i\theta} \colon e^{i\theta}\in I, \, 1-|I| < r < 1 \},\]
    where $|I|$ denotes the arc length of $I$. Then, $\nu$ is Carleson for $H^2(\disk)$ if and only if there exists a constant $C>0$ such that 
    \begin{equation}
        \label{carmeasH2eq}
        \nu(S(I))\leq C|I|, \qquad I\subset\T.
    \end{equation}
\end{propo}

Carleson measures of $D_\zeta$ have been characterized in \cite{20} in terms of Carleson measures of $H^2(\disk)$, as follows: 
\begin{propo}
\label{carmeasDz}
    Let $\nu$ be a finite positive Borel measure on $\disk$. Then, $\nu$ is a Carleson measure for $\D_\z$ if and only if the measure $|z-\z|^2 \de \nu(z)$ is Carleson for $H^2(\disk)$.
\end{propo}
Note that every Carleson measure of $\D_\z$ has to be finite, since $1\in\D_\z$. Having mentioned these preliminary facts, we can state our result. 

\begin{cor}
\label{Carmeasureinc}
    Let $u$ be an inner function with $\sigma(u) \neq \T$, and $\nu$ a finite positive Borel measure on $\disk$. If there exists $\zeta\in\T\setminus\sigma(u)$ such that $|z-\z|^2\de\nu(z)$ is a Carleson measure for $H^2(\disk)$, then $\nu$ is a Carleson measure for the model space $K_u$.
\end{cor}

\begin{proof}
Since $\z\notin\sigma(u)$, by Theorem \ref{KuinDz} the embedding $K_u\hookrightarrow\D_\z$ holds. Also, by Proposition \ref{carmeasDz}, the measure $\nu$ is a Carleson measure for $\D_\z$. Then, for every $f\in K_u$ it holds
\[\int_\disk |f|^2 \de\nu \leq C \|f\|_{\D_\z}^2 \leq C' \|f\|_{K_u}^2,\]
for some positive constants $C,C'$, meaning that $\nu$ is a Carleson measure for $K_u$.
\end{proof}

We conclude this part with an example of a Carleson measure for $\D_1$ (and thus for every model space $K_u$ with $1\notin\sigma(u)$) which is not Carleson for $H^2(\disk)$.

\begin{example}
Let $\nu$ be the measure defined on Borel sets of $\disk$ as
\[ \nu(A) := \int_{A\cap[0,1]} \frac{1}{\sqrt{1-s}}\de s.\]
We use the characterization in Proposition \ref{carmeasH2} to prove that $\nu$ is not a Carleson measure for $H^2(\disk)$. For $\delta>0$, consider the arc $I_\delta$ centered at $1$ with arc length $\delta$. One can compute the measure of the Carleson boxes $S(I_\delta)$ and obtain
\[\nu\big(S(I_\delta)\big) = \int_{1-\delta}^1 \frac{1}{\sqrt{1-s}}\de s = 2\sqrt{\delta}, \]
showing that the bound in \eqref{carmeasH2eq} cannot hold as $\delta \to 0$. However, the measure $\nu$ is a Carleson measure for the local Dirichlet space $\D_1$. We use Proposition \ref{carmeasDz}, and because of the definition of $\nu$ it suffices to consider only the arcs that contain $1$, and one can show that the measure $|z-1|^2\de\nu(z)$ satisfies \eqref{carmeasH2eq}.
\end{example}

  We move now to the description of multipliers.
\begin{defi}
    Let $H_1,H_2$ be Hilbert spaces of holomorphic functions on $\disk$. The multipliers from $H_1$ to $H_2$ are defined as
    \[ 
    M(H_1,H_2) := \{\phi\in\Hol(\disk) \colon \phi H_1 \subseteq H_2\}.
    \]
    When $H_1=H_2$ we simply write $M(H_1)$.
\end{defi}

The multiplier algebra $M(\D_\z)$ of the local Dirichlet space is characterized as follows. This result follows from Proposition $3.1$ of \cite{fricain2018multipliers}. For sake of completeness, we provide an explicit proof.

\begin{lem}
\label{multDz}
For $\z\in\T$, the multiplier algebra of $\D_\z$ is $\D_\z \cap H^\infty(\disk)$.
\end{lem}

\proof
The fact that the multipliers of $\D_\z$ are in $\D_\z \cap H^\infty(\disk)$ follows from the standard argument which holds for many other reproducing kernel Hilbert spaces of analytic functions, see for example Proposition $3.1$ in \cite{modelmultipliers}. Let us move to the other inclusion: let $\phi \in \D_\z \cap H^\infty(\disk)$, and let $f\in\D_\z$. In light of the characterization in \eqref{second defn}, there exist functions $\eta, g\in H^2(\disk)$ such that 
\begin{equation}
\label{multDzeq}
\phi(z) = \phi(\z)+(z-\z)\eta(z), \qquad f(z)=f(\z) + (z-\z)g(z), \qquad z\in\disk.
\end{equation}
Then, for $z\in\disk$ it holds 
\begin{align*}
\phi(z)f(z) &= \big( \phi(\z)+(z-\z)\eta(z) \big) \big( f(\z) + (z-\z)g(z) \big) \\
&= \phi(\z)f(\z) + (z-\z)[ \phi(\z)g(z)+\eta(z)f(\z) + (z-\z) \eta(z)g(z)]. 
\end{align*}
Again by \eqref{second defn}, membership of the product $\phi f$ in $\D_\z$ is equivalent to the membership in $H^2(\disk)$ of the function
\[ \phi(\z)g(z)+\eta(z)f(\z) + (z-\z) \eta(z)g(z).\] Since $\eta,g\in H^2(\disk)$, it suffices to show that $(z-\z) \eta(z)g(z)$ belongs to $H^2(\disk)$, and this follows from \eqref{multDzeq} and the assumption that $\phi\in H^\infty(\disk)$, for
\[ (z-\z) \eta(z)g(z) = \big(\phi(z)-\phi(\z)\big)g(z). 
\qed\]

In \cite{modelmultipliers}, multipliers between model spaces are studied. It is shown that $M(K_u)=\C$, meaning that every function multiplying any model space into itself must be constant. Furthermore, multipliers from model spaces to the Hardy space $H^2(\disk)$ are characterized in terms of a Carleson condition on the unit circle. More precisely, $\phi\in\M\big(K_u,H^2(\disk)\big)$ if and only if the measure $|\phi|^2\!\de m$ is a Carleson measure for $K_u$, i.e. there exists a constant $C>0$ such that
\[ \int_\T |f\phi|^2 \de m \leq C \|f\|_{K_u}^2, \qquad f\in K_u .\] Assuming the inclusion $K_u\subseteq \D_\z$, the local Dirichlet space $\D_\z$ is an intermediate space between $K_u$ and $H^2(\disk)$. This is reflected in our following multiplier theorem. 
\begin{thm}
\label{thmmoltKuDz}
    Let $u$ be an inner function, $\z\in\T$ such that $\z\notin\sigma(u)$, and $\phi\in\emph{Hol}(\disk)$. Then $\phi$ is a multiplier from $K_u$ to $\D_\z$ if and only if the measure $|\phi|^2\!\de m$ is Carleson for $K_u$ and $\phi$ belongs to $\D_\z$.
\end{thm}

\begin{proof}[Proof of Theorem \ref{thmmoltKuDz} ]
Let us assume that $\phi\in M(K_u,\D_\z)$. Then, in particular, the measure $|\phi|^2 \de m$ is a Carleson measure for $K_u$, so it suffices to show that every multiplier from $K_u$ to $\D_\z$ belongs to $\D_\z$. If $u(0)=0$, then $1\in K_u$, implying that the multiplier $\phi$ belongs to $\D_\z$. If $u(0)\neq 0$, we consider the kernel 
\[ 
    k_0^u= 1-\overline{u(0)}u.
    \]
    Using Theorem \ref{localdouglasthm}, one can check that $1/k_0^u \in H^\infty(\disk)\cap\D_\z$, so that by Lemma \ref{multDz} the function $1/k_0^u$ is a multiplier of $\D_\z$. Thus, 
    \[ 
    \phi = \frac{1}{k_0^u} \, \phi k_0^u \in \D_\z
    \]
    which implies the statement.
    Let us now prove the other implication.
    We assume that $|\phi|^2 \de m$ is a Carleson measure for $K_u$ and that $\phi$ belongs to $\D_\z$. Since $\phi\in M\big(K_u,H^2(\disk)\big)$, for every $f\in K_u$ the product $\phi f$ belongs to $H^2(\disk)$. We compute the local Dirichlet integral.
    \begin{align*}
        D_\z(f\phi) &= \int_\T \bigg| \frac{f(\la)\phi(\la)-f(\z)\phi(\z)}{\la-\z}\bigg|^2 \de m(\la) \\
        &= \int_\T \bigg| \frac{f(\la)\phi(\la)-\phi(\la)f(\z)+\phi(\la)f(\z) - f(\z)\phi(\z)}{\la-\z}\bigg|^2 \de m(\la) \\
        &\leq \int_\T |\phi(\la)|^2 \bigg| \frac{f(\la)- f(\z)}{\la-\z} \bigg|^2 \de m(\la) + |f(\z)|^2 \int_\T \bigg| \frac{\phi(\la)- \phi(\z)}{\la-\z} \bigg|^2 \de m(\la) \\
        &= \bigg\|\frac{f- f(\z)}{\cdot-\z} \bigg\|_{L^2(|\phi|^2\!\de m)}^2 + |f(\z)|^2\bigg\| \frac{\phi- \phi(\z)}{\cdot-\z} \bigg\|_{H^2}^2 \\
        &\leq C \bigg\|\frac{f- f(\z)}{\cdot-\z} \bigg\|_{K_u}^2 + |f(\z)|^2 D_\z(\phi) \\
        &\leq \big( C+D_\z(\phi)\big) \|f\|_{\D_\z}^2, 
    \end{align*}
concluding the proof.
\end{proof}
    
It is natural to ask whether the condition in Theorem \ref{thmmoltKuDz} guarantees the boundedness of the multipliers, in other words, whether $M(K_u,\D_\z)$ is contained or not in $H^\infty(\disk)$.  The answer to this question is negative. Considering the simplest case $u(z)=z$, one has that $K_u = \mathbb{C}$, and therefore $M(K_u,\D_\zeta)=\D_\zeta$, which contains unbounded functions. 

\vspace{11 pt} 
\section{The embedding $K_u \hookrightarrow \Dm$}
In this section we study the embedding of $K_u$ into $\Dm$, for an arbitrary measure $\mu$. In this case, the sufficient condition we obtain is different from the necessary one. We start with the proof of Theorem \ref{KuinDm}.

\begin{proof}[Proof of Theorem \ref{KuinDm}.]
By assumption, $supp(\mu)$ and $\sigma(u)$ are disjoint compact sets, therefore
\[ 
\delta := \text{dist}\big(\supp(\mu),\sigma(u)\big) >0.
\]
We consider the open set 
\[ U:=\bigcup_{x\in\sigma(u)} \bigg\{z\in\C\colon |z-x|<\frac{\delta}{2} \bigg\}.
\]
We split the harmonically weighted Dirichlet integral into
\[ 
D_\mu(f) = \frac{1}{\pi} \int_{\disk \cap U} |f'|^2 P\mu\de A+ \frac{1}{\pi} \int_{\disk \setminus U} |f'|^2 P\mu\de A.
\]
For the first summand, we use a classical Littlewood-Paley estimate, see Proposition $3.2$ in \cite{garnett}:
\begin{align*}
\frac{1}{\pi}\int_{\disk \cap U} |f'|^2 P\mu\de A &=\frac{1}{\pi}\int_{\disk \cap U} |f'(z)|^2  (1-|z|^2)  \bigg( \int_{\supp(\mu)} \frac{\de \mu(\z)}{|z-\z|^2}\bigg) \de A(z) \\
&\leq \frac{4}{\pi \delta^2}\mu(\T) \int_{\disk} |f'(z)|^2  (1-|z|^2) \de A(z) \\
&\leq \frac{2}{\delta^2}\mu(\T)\int_\T |f(\la)-f(0)|^2\de m(\la) \\
&\leq \frac{8}{\delta^2}\mu(\T) \|f\|_{H^2}^2.
\end{align*}
For the second summand, we recall that every function in the model space $K_u$ admits an analytic extension across $\T\setminus\sigma(u)$. Hence, we have
\begin{align*}
\frac{1}{\pi} \int_{\disk \setminus U} |f'(z)|^2 P\mu(z)\de A(z) &\leq \max_{\overline{\disk} \setminus U} |f'|\, \frac{1}{\pi}\int_\disk \int_\T \frac{1-|z|^2}{|z-\z|^2} \de \mu(\z) \de A(z) \\
&= \max_{\overline{\disk} \setminus U} |f'| \, \int_\T \de \mu(\z) \\
&= \max_{\overline{\disk} \setminus U} |f'|\,  \mu(\T).
\end{align*}
We have proved that for every $f\in K_u$
\[  D_\mu(f) \leq \frac{8}{\delta^2}\mu(\T) \|f\|_{H^2}^2 + \max_{\overline{\disk} \setminus U} |f'| \mu(\T) <\infty.  
\]
The boundedness of the embedding $K_u \hookrightarrow\Dm$ follows from the closed graph theorem.
\end{proof}

Now we prove Theorem \ref{KuinDmnec}, giving a necessary condition for the considered embedding.

\begin{proof}
For the proof, we introduce the function $V_\mu\colon\C\to[0,+\infty]$ defined as
\begin{equation*}
V_\mu(\om):= \int_{\T}\frac{1}{|\z-\om|^2}\de\mu(\z), \qquad \om\in\C.
\end{equation*}
First, we prove that $V_\mu$ is bounded on the boundary spectrum $\sigma(u)$, which we can assume to be non-empty without loss of generality. Let $C>0$ be a constant such that
\[D_\mu(f) \leq C\|f\|_{H^2}^2, \qquad f\in K_u.\]
Let $\la\in\sigma(u)$ and, as we did in the proof of Theorem \ref{KuinDzb}, let us consider a sequence $(\om_n)_n$ in $\disk$ such that $u(\omega_n) \to 0$ as $\omega_n \to \lambda$. 
By the disintegration formula in \eqref{eqfubini} and the lower estimate for $\D_\z(k_n)$ obtained in the proof of Theorem \ref{KuinDzb}, we have that
\begin{align*}
    C\|k_n\|_{H^2}^2&\geq    D_\mu(k_n)=  \int_{\T} D_\z(k_n) \de \mu(\z) \\
    &\geq   \int_{\T} \|k_n\|_{H^2}^2 \frac{\big(1-|u(\om_n)|\big)^2 \big(|\om_n| - |u(\om_n)|\big)^2}{|\z-\om_n|^2 \big(1-|u(\om_n)|^2\big)} \de\mu(\z) \\
    &=   \|k_n\|^2_{H^2} \frac{\big(1-|u(\om_n)|\big)\big(|\om_n| - |u(\om_n)|\big)^2}{1+|u(\om_n)|} \int_{\T} \frac{1}{|\zeta -\omega_n|^2}  \de\mu(\z) .
\end{align*}
Hence, by Fatou's Lemma, it holds that 
\[ C \geq \liminf_{n} \int_{\T} \frac{1}{|\zeta -\omega_n|^2}  \de\mu(\z) \geq \int_{\T} \frac{1}{|\zeta -\la|^2}  \de\mu(\z)=V_\mu(\la),\]
which proves that $\sup_{\la \in \sigma(u)} V_\mu(\la) <\infty$. Now the theorem follows from the fact that $V_\mu = \infty$ $\mu$-a.e. on $\T$ and therefore, necessarily, we have that $\mu\big(\sigma(u)\big)=0$.
\end{proof}

\begin{oss}
    We note that a similar necessary condition holds also for the embedding $H(b)\hookrightarrow \Dm$. Let $b_i$ be the inner factor associated to the bounded function $b$, and we consider a point $\z\in\sigma(b_i)$. If $\lim_n \om_n = \z$ and $\lim_n |b(\om_n)|=0$ we note that
    \[
    C\|k_n\|_{b}^2\geq \int_\T D_\zeta(k_n)d\mu(\zeta)\geq \|k_n\|^2_b \frac{\big(1-|b(\om_n)|\big)\big(|\om_n| - |b(\om_n)|\big)^2}{1+|b(\om_n)|} \int_{\T} \frac{1}{|\zeta -\omega_n|^2}  d\mu(\z),
    \]
    and once again by Fatou's Lemma we conclude that $V_\mu$ is bounded on $\sigma(b_i)$ and therefore $\mu\big(\sigma(b_i)\big)=0$.
\end{oss}

  We conclude this section discussing the compactness of the embeddings. Due to the trivial norm inequality $\|\cdot\|_{K_u}\leq\|\cdot\|_\mu$, the compactness of the embedding $K_u\hookrightarrow\Dm$ implies the compactness of the identity map $I_{K_u}$. Therefore, it is easy to see that the embedding $K_u\hookrightarrow\Dm$ is compact if and only if $K_u$ is finite dimensional, that is, if and only if $u$ is a finite Blaschke product.

\section{Final remarks and open questions}
Given $\mu$ a finite positive Borel measure on $\T$ and an inner function $u$, we have provided a sufficient condition and a necessary condition for the embedding $K_u\hookrightarrow\Dm$, respectively $\supp(\mu)\cap\sigma(u)=\emptyset$ and $\mu\big(\sigma(u)\big)=0$. If $\mu=\delta_\z$, both these conditions are equivalent to $\z\notin\sigma(u)$. For the Lebesgue measure, the two conditions do not coincide, but the sufficient one is also necessary. This is because, if the inclusion $K_u \hookrightarrow \D=\D(m)$ holds, then necessarily $u$ belongs to $\D$: taking  $\omega\in\mathbb{D}$ such that $u(\omega)\neq 0$, one has
\[u(z)=\frac{1}{\overline{u(\omega)}}\left[1-\left(1-\overline{\omega}z\right) k_\omega^u(z)\right],\qquad z\in\mathbb{D},\]
so that $u = \overline{u(\omega)}^{-1} \left(1-\left(I-\overline{\omega}S\right)k_\omega^u\right) \in \D$.
However, it is shown in \cite{primer} that the only inner functions in the classical Dirichlet space are finite Blaschke products, resulting in the boundary spectrum $\sigma(u)$ being empty. In future works we will investigate whether the sufficient condition $\supp(\mu)\cap\sigma(u)=\emptyset$ is in general necessary as well for the embedding $K_u\hookrightarrow\Dm$. For the time being, we leave this as an open problem. 

\section*{Acknowledgements} 
The authors want to thank Professor Maria Nowak for the interest in discussing some problems addressed in this article. They also want to thank Professor Marco M. Peloso for his constant help and guidance.


\bibliographystyle{plain}
\bibliography{mybibliography}

\end{document}